\documentclass[12pt]{article}
\usepackage{latexsym, amssymb, amsmath, amscd, amsfonts, epsfig, graphicx, colordvi, ifpdf}

\newtheorem{thm}{Theorem}[section]

\newtheorem{defi}[thm]{Definition}
\newtheorem{lem}[thm]{Lemma}
\newtheorem{core}[thm]{Corollary}
\newtheorem{pro}[thm]{Proposition}

\def\pf{\noindent{\it Proof.} }
\def\qed{\nopagebreak\hfill{\rule{4pt}{7pt}}
\medbreak}

\def\qed{\nopagebreak\hfill{\rule{4pt}{7pt}}
\medbreak}

\newenvironment{kst}
{\setlength{\leftmargini}{2.4\parindent}
\begin{itemize}
\setlength{\itemsep}{-0.5mm}} {\end{itemize}}

\title{The combinatorics of $k$-marked Durfee symbols}
\author{ Kathy Qing Ji \\
 \small Center for Combinatorics, LPMC-TJKLC
\\[-0.8ex]
 \small Nankai University, Tianjin 300071, P.R. China\\
\small \texttt{ji@nankai.edu.cn} }

\begin{document}
\maketitle

\vskip 6mm \noindent {\bf Abstract.}  Andrews recently introduced
$k$-marked Durfee symbols which are connected to moments of
Dyson's rank. By these connections, Andrews  deduced  their
generating functions and some combinatorial properties and left
their purely combinatorial proofs as open problems. The primary
goal of this article is to provide combinatorial proofs in answer
to Andrews' request. We obtain a relation between $k$-marked
Durfee symbols and Durfee symbols by constructing bijections, and
all identities on $k$-marked Durfee symbols given by Andrews could
follow from this relation. In a similar manner, we also prove the
identities due to Andrews on $k$-marked odd Durfee symbols
combinatorially, which resemble ordinary $k$-marked Durfee symbols
with a modified subscript and with odd numbers as entries.



\section{Introduction}

We will adopt the terminology on partitions in Andrews
\cite{Andrews-1976}.  {\it A partition} $\lambda$ of a positive
integer $n$ is a finite nonincreasing sequence of positive
integers $(\lambda_1,\,\lambda_2,\ldots,\,\lambda_r)$ such that
$\sum_{i=1}^r\lambda_i=n.$ Then $\lambda_i$ are called the parts
of $\lambda$, $\lambda_1$ is its largest part. The number of parts
of  $\lambda$ is called the length of $\lambda$, denoted by
$l(\lambda).$ The weight of $\lambda$ is the sum of parts of
$\lambda$, denoted by $|\lambda|.$

{\it The rank of a partition} $\lambda$ introduced  by Dyson
\cite{Dyson-1944} is defined as the largest part minus the number
of parts, which is usually denoted by
$r(\lambda)=\lambda_1-l(\lambda).$ Let $N(m;n)$ denote the number
of partitions of $n$ with rank $m$, we have

\begin{thm}[Dyson]\label{gf-r}
The generating function for $N(m;n)$ is given by
\begin{equation}
\sum_{n=0}^{+\infty}
N(m;n)q^n=\frac{1}{(q;q)_{\infty}}\sum_{n=1}^{+\infty}(-1)^{n-1}q^{n(3n-1)/2+|m|n}(1-q^n),\quad
|q|<1.
\end{equation}
\end{thm}
where $(a;q)_n=\prod_{j=0}^{n-1}(1-aq^j)$ and
$(a;q)_{\infty}=\lim_{n\rightarrow\infty}(a;q)_n$. The identity
\eqref{gf-r} was first discovered by Dyson \cite{Dyson-1944} in
1944 and first proved by Atkin and Swinnerton-Dyer
\cite{Atkin-Swinnerton-Dyer-1954}. Later Dyson \cite{Dyson-1969}
gave a simple combinatorial argument of it. We refer to
\cite[p.63]{Berkovich-Garvan-2002} for more details.
\begin{defi} For a nonnegative integer $n$, a Durfee symbol of $n$  is a
 two-row array with a subscript
\begin{equation}
(\alpha,\beta)_D=\left(\begin{array}{cccc}
\alpha_1&\alpha_2&\ldots&\alpha_s\\
\beta_1&\beta_2&\ldots&\beta_t
\end{array} \right)_{D}
\end{equation}
where $D\geq \alpha_1\geq \alpha_2\geq \cdots \geq \alpha_s>0$,
$D\geq \beta_1 \geq \beta_2 \geq \cdots \geq \beta_t>0$ and
$n=\sum_{i=1}^s \alpha_i+\sum_{i=1}^t \beta_i+D^2$.
\end{defi}
For example, there are $5$ Durfee symbols of $4$.
\[\left(\begin{array}{ccc}
1&1&1\\
&&
\end{array}
\right)_1\ \left(\begin{array}{cc}
1&1\\
1&
\end{array}
\right)_1\ \left(\begin{array}{cc}
1&\\
1&1
\end{array}
\right)_1\ \left(\begin{array}{ccc}
\\
1&1&1
\end{array}
\right)_1\ \left(\begin{array}{c}
\ \\
\
\end{array}
\right)_2\] The difference of the lengths of $\alpha$ and $\beta$
is called {\it rank} of Durfee symbols $(\alpha,\,\beta)_D$. We
use $\mathcal{D}_1(m;n)$ to denote the number of Durfee symbols of
$n$ with rank $m$. Andrews \cite[Section 3]{Andrews-07-a} showed
that by constructing a bijection
\begin{thm}[Andrews]\label{dy-op}
The number of ordinary partitions of $n$ with rank equal to $m$ is
equal to the number of Durfee symbols of $n$ with rank equal to
$m$, that is
\begin{equation}
N(m;n)=\mathcal{D}_1(m;n).
\end{equation}
\end{thm}
Andrews \cite{Andrews-07-a} introduced the $k$th symmetrized rank
moment $\eta_k(n)$, defined by
\begin{equation}\label{sym-ran-mom}
\eta_k(n)=\sum_{m=-\infty}^{{+\infty}}
{m+\lfloor\frac{k-1}{2}\rfloor \choose k}N(m,n),
\end{equation}
which are linear combinations of the $k$th rank moments $N_k(n)$
\begin{equation}
N_k(n)=\sum_{m=-\infty}^{{+\infty}} m^kN(m,n),
\end{equation}
 considered by Atkin and Garvan \cite{Atkin-Garvan-03}.  To give a
combinatorial explanation of \eqref{sym-ran-mom}, Andrews
\cite{Andrews-07-a} introduced $k$-marked Durfee symbols, which
can be thought of as the generalized Durfee symbols.

\begin{defi}\label{de-MDS}
A $k$-marked Durfee symbol of $n$ is composed of $k$ pairs of
partitions and a subscript,  which is defined as
 \[\eta=\left(\begin{array}{cccc}
\alpha^k,&\alpha^{k-1},&\ldots,&\alpha^1\\[2pt]
\beta^k,&\beta^{k-1},&\ldots,&\beta^1
\end{array}
\right)_D,\] where $\alpha^i$ {\rm(}resp. $\beta^i${\rm)}
represents a  partition and
$\sum_{i=1}^k(|\alpha^i|+|\beta^i|)+D^2=n$. Furthermore, the
partitions $\alpha^i$ and $\beta^i$ must satisfy the following
three conditions where  $\alpha^i_1$ {\rm(}resp. $\beta^i_1${\rm)}
is the largest part of the partition $\alpha^i$ {\rm(}resp.
$\beta^i${\rm)} and $\alpha^i_{l(\alpha^i)}$ {\rm(}resp.
$\beta^i_{l(\beta^i)}${\rm)} is the smallest part of the partition
$\alpha^i$ {\rm(}resp. $\beta^i${\rm)}.
\begin{itemize}
\item[{\rm (1)}]For $1\leq i<k$, $\alpha^i$ must be non-empty
partition, while $\alpha^k$ and  $\beta^i$ could be empty{\rm ;}

\item[{\rm (2)}] $\beta^{i-1}_1\leq \alpha^{i-1}_1 \leq
\beta^{i}_{l(\beta^{i})}$ for $2\leq i\leq k${\rm ;}

\item[{\rm (3)}]$\beta^k_1,\alpha^k_1 \leq D.$
\end{itemize}
\end{defi}
Clearly,  $1$-marked Durfee symbol is just Durfee symbol.

 Let $\eta=\left(\begin{array}{cccc}
\alpha^k,&\alpha^{k-1},&\ldots,&\alpha^1\\[2pt]
\beta^k,&\beta^{k-1},&\ldots,&\beta^1
\end{array}
\right)_D$ be a $k$-marked Durfee symbol. The pair of partitions
$\left(\begin{array}{l}
\alpha^i\\
\beta^i
\end{array}
\right)$ is called the $i$th vector of $\eta$. We define
$\rho_i(\eta)$, the $i$th rank of $\eta$ by
\[\rho_i(\eta)=\left\{
\begin{array}{ll}
l(\alpha^i)-l(\beta^i)-1 \ \ &\text{ for }\ 1\leq i<k,\\[5pt]
l(\alpha^k)-l(\beta^k) \ \ &\text{ for }\ i=k.
\end{array}\right.
\]
For example, $\left(\begin{array}{ccccccc}
4_3&4_3&3_2&3_2&2_2&2_1\\
&5_3&&3_2&2_2&2_1
\end{array}\right)_5$ is a $3$-marked Durfee symbol of $55$ where
$\alpha^3=(4, 4),\,\alpha^2=(3,3,2),\alpha^1=(2),$ and
$\beta^3=(5),\beta^2=(3,2),\beta^1=(2).$ The first rank is $-1$,
the second rank is $0$, and the third rank is $1.$

Let $\mathcal{D}_k(m_1,m_2,\ldots,m_k;n)$ denote the number of
$k$-marked Durfee symbols of $n$ with $i$th rank equal to $m_i$
and $\mathcal{D}_k(n)$ denote the number of $k$-marked Durfee
symbols of $n$, it's clear to see that
\begin{equation}
\mathcal{D}_k(n)=\sum_{m_1,\ldots,m_{k}=-\infty}^{+\infty}
\mathcal{D}_k(m_1,m_2,\ldots,m_k;n).
\end{equation}

In his recent work \cite{Andrews-07-a}, Andrews used the
connections between $k$-marked Durfee symbols and the symmetrized
rank moments \eqref{sym-ran-mom} to find identities relating the
generating function for $\mathcal{D}_k(m_1,m_2,\ldots,m_k;n)$, as
well as to deduce some combinatorial properties and Ramanujan-type
congruences for $k$-marked Durfee symbols. At the end of the
paper, Andrews proposed  a variety of serious questions which fall
into 3 basic groups: combinatorial, asymptotic and congruential.
The recent works \cite{Bringmann-1, Bringmann-2} by Kathrin
Bringmann, Frank Garvan, and Karl Mahlburg  focused on the study
relating asymptotical and congruential properties of $k$-marked
Durfee symbols.  They used the automorphic properties to prove the
existence of infinitely many congruences for $k$-marked Durfee
symbols. The primary goal of this article is to answer Andrews'
request on combinatorics (Problems 1-4 and 6-9 on page 39 of
\cite{Andrews-07-a}). We will give combinatorial proofs of the
identities relating the generating function for
$\mathcal{D}_k(m_1,m_2,\ldots,m_k;n)$ and combinatorial properties
for $k$-marked Durfee symbols. To be specific, we first derive the
following partition identity \eqref{main-e} by constructing
bijections, which gives a relation between $k$-marked Durfee
symbols and Durfee symbols. We then show that all identities on
$k$-marked Durfee symbols given  by Andrews ( \cite[Problems
1-4]{Andrews-07-a}) could follow from this identity. We then use
the similar method to study the identities of Andrews on
$k$-marked odd Durfee symbols ( \cite[Problems
6-9]{Andrews-07-a}), which resemble ordinary Durfee symbols with a
modified subscript and with odd numbers as entries.

\begin{thm}\label{main}
For $k\geq 2$, we have
\begin{equation}\label{main-e}
\mathcal{D}_k(m_1,m_2,\ldots,m_k;n)=\sum_{j=0}^{+\infty}{j+k-2
\choose k-2}N\left(\sum_{i=1}^k|m_i|+2j+k-1;n\right).
\end{equation}
\end{thm}
The paper is organized as follows. In Section 2, we consider the
relation between $k$-marked Durfee symbols and Durfee symbols and
prove Theorem \ref{main}. To this end, we introduce a special
class of $k$-marked Durfee symbols, which we call  $k$-marked
strict shifted Durfee symbols since each of their vectors except
for the $k$th vector is a two-line strict shifted plane partition.
We deduce the desired relation by building the connections between
$k$-marked strict shifted Durfee symbols, Durfee symbols, and
$k$-marked Durfee symbols respectively. In Section 3, we give
combinatorial proofs of the identities due to Andrews on
$k$-marked Durfee symbols with the help of Theorem \ref{main}. In
particular, the symmetry of $k$-marked Durfee symbols
(\cite[Corollary 12]{Andrews-07-a}) could be thought of as a
direct consequence  of Theorem \ref{main}. Section 4 is devoted to
the study of  $k$-marked odd Durfee symbols.

\section{$k$-marked strict shifted Durfee symbols} In this section, we will
establish  the relations between $k$-marked strict shifted Durfee
symbols,  Durfee symbols, and $k$-marked Durfee symbols
respectively, and then deduce Theorem \ref{main}. We begin by
defining $k$-marked strict shifted Durfee symbols.

To define these objects, we need recall the concept of strict
shifted plane partitions introduced by Andrews
\cite{Andrews-1979}. Mills, Robbins and Rumsey
\cite{Mills-Robbins-Rumsey} showed that strict shifted plane
partitions whose row lengths are equal to row leaders are
bijective to cyclically symmetric plane partitions.

A strict shifted plane partition of $n$ is an array $\pi=(a_{ij})$
of positive integers defined only for $j\geq i$, which has
non-increasing rows and strictly decreasing columns, such that the
sum of its elements is $n$. Such an array can be written
\[\begin{array}{ccccclll}
a_{11}&a_{12}&a_{13}&a_{14}&\cdots&&&a_{1\mu_1}\\
&a_{22}&a_{23}&a_{24}&\cdots&&a_{2\mu_2}\\
&&&&\cdots\cdots&& \\
&&&a_{rr}&\cdots&a_{r\mu_r}
\end{array}
\]
where $a_{ij}\geq a_{i,j+1}$ and $a_{ij}> a_{i+1,j}$, $\mu_1\geq
\mu_2\geq\cdots\geq \mu_r$ and $\sum_{i,j}a_{ij}=n$.

We can regard a two-lined strict shifted partition of $n$ as a
pair of partitions $(\alpha,\beta)$ of $n$, where
$\alpha=(\alpha_1,\alpha_2,\cdots,\alpha_r),\,\beta=(\beta_1,\beta_2,\cdots,\beta_s),r>
s,$ and $\alpha_{i+1}>\beta_{i}$
 for $i=1,2,\ldots,s$. For example $\left(\begin{array}{cccccc}
3&3&3&2&2&1\\
&2&1&1&1
\end{array}
\right)$ is a two-lined strict shifted partition.

 A $k$-marked Durfee symbol $\eta=\left(\begin{array}{cccc}
\alpha^k,&\alpha^{k-1},&\ldots,&\alpha^1\\[2pt]
\beta^k,&\beta^{k-1},&\ldots,&\beta^1
\end{array}
\right)_D$ is said to be {\it strict shifted} if all of its
vectors except for the $k$th vector $(\alpha^k,\beta^k)$ are
two-lined strict shifted partitions. For example, $3$-marked
Durfee symbol in \eqref{example} is strict shifted. Let
$\mathcal{D}^{ss}_k(m_1,m_2,\ldots,m_k;n)$ denote the number of
$k$-marked strict shifted Durfee symbols of $n$ with $i$th rank
equal to $m_i$.

 We now build a connection between $k$-marked strict shifted Durfee
symbols of $n$ and Durfee symbols of $n$.
\begin{thm}\label{r-kms-ds}
Given $k$ nonnegative integers $m_1,\,m_2,\ldots,m_{k}$, there is
a bijection $\Phi$ between the set of $k$-marked strict shifted
Durfee symbols of $n$ with $i$th rank equal to $m_i$ and the set
of Durfee symbols of $n$ with rank equal to $\sum_{i=1}^km_i+k-1$.
\end{thm}
\pf {\bf The map $\Phi$:} Let $\eta=\left(\begin{array}{cccc}
\alpha^k,&\alpha^{k-1},&\ldots,&\alpha^1\\[2pt]
\beta^k,&\beta^{k-1},&\ldots,&\beta^1
\end{array}
\right)_D$ counted by $\mathcal{D}^{ss}_k(m_1,\ldots,m_k;n)$, we
then obtain a Durfee symbol $(\gamma,\delta)_{D}$ when remove all
subscripts of $\eta$. Obviously, the resulting Durfee symbol
$(\gamma,\delta)_{D}$ is enumerated by $
\mathcal{D}_1(\sum_{i=1}^km_i+k-1;n).$

{\bf The reverse map $\Phi^{-1}$:} Let $(\gamma,\delta)_{D}$ be
counted by $\mathcal{D}_1(\sum_{i=1}^km_i+k-1;n)$, we will
construct a $k$-marked Durfee symbol
$\eta=\left(\begin{array}{cccc}
\alpha^k,&\alpha^{k-1},&\ldots,&\alpha^1\\[2pt]
\beta^k,&\beta^{k-1},&\ldots,&\beta^1
\end{array}
\right)_{D'}$ whose $i$th rank equal to $m_i$. Let
\[\left(\begin{array}{l}
\gamma\\
\delta
\end{array}\right)_D=\left(\begin{array}{cccc}
\gamma_1&\gamma_2&\ldots&\gamma_l\\[2pt]
\delta_1&\delta_2&\ldots&\delta_s
\end{array}
\right)_D,\] where $D\geq \gamma_1\geq \gamma_2\geq \cdots\geq
\gamma_l$ and $D\geq \delta_1\geq \delta_2\geq \cdots \geq
\delta_s$, we assume that $\delta_{j}=0$ for $j\geq s+1$. Note
that $l-s=\sum_{i=1}^km_i+k-1$.

We now split $(\gamma, \delta)$ to generate  the $k$th vector
$(\alpha^k,\beta^k)$ of $\eta$.
 Let $j$ be largest nonnegative integer such that
$\delta_{j}\geq \gamma_{m_k+j+1}$,  that is for any $i\geq j+1$,
we have $\delta_{i}<\gamma_{m_k+i+1}$. Let
\[\left(\begin{array}{l}
\alpha^k\\
\beta^k
\end{array}\right)=\left(\begin{array}{cccc}
\gamma_1&\gamma_2&\ldots&\gamma_{m_k+j}\\[2pt]
\delta_1&\delta_2&\ldots&\delta_{j}
\end{array}
\right),\] and
\[\left(\begin{array}{l}
\gamma'\\
\delta'
\end{array}\right)=\left(\begin{array}{cccc}
\gamma'_1&\gamma'_2&\ldots&\gamma'_{l'}\\[2pt]
\delta'_1&\delta'_2&\ldots&\delta'_{s'}
\end{array}
\right),\] where $\gamma'_i=\gamma_{m_k+j+i}$,
$\delta'_i=\delta_{j+i}$ for $i\geq 1$. Obviously,
$l(\alpha^k)-l(\beta^k)=m_k$. Furthermore,  $(\gamma',\delta')$ is
a strict shifted partition from the fact that for any $i\geq 1$,
$\delta_{i+j}<\gamma_{m_k+i+j+1}$ and
$l'-s'=\sum_{i=1}^{k-1}m_i+k-1.$

We continue to split $(\gamma', \delta')$ to construct the
$(k-1)$th vector $(\alpha^{k-1},\beta^{k-1})$  of $\eta$. Let $j$
be largest nonnegative integer such that $\delta'_{j}\geq
\gamma'_{m_{k-1}+j+2}$, we then let
\[\left(\begin{array}{l}
\alpha^{k-1}\\
\beta^{k-1}
\end{array}\right)=\left(\begin{array}{cccc}
\gamma'_1&\gamma'_2&\ldots&\gamma'_{m_{k-1}+j+1}\\[2pt]
\delta'_1&\delta'_2&\ldots&\delta'_{j}
\end{array}
\right),\]

and
\[\left(\begin{array}{l}
\gamma''\\
\delta''
\end{array}\right)=\left(\begin{array}{cccc}
\gamma''_1&\gamma''_2&\ldots&\gamma''_{l''}\\[2pt]
\delta''_1&\delta''_2&\ldots&\delta''_{s''}
\end{array}
\right),\] where $\gamma''_i=\gamma'_{m_{k-1}+j+i+1}$,
$\delta''_i=\delta'_{j+i}$ for $i\geq 1$. Clearly,
$l(\alpha^{k-1})-l(\beta^{k-1})=m_{k-1}+1$ and
$(\alpha^{k-1},\beta^{k-1})$ is a strict shifted partition for
$(\gamma',\delta')$ is strict shifted. Observe that
$\delta'_{j+i}< \gamma'_{m_{k-1}+i+j+2}$ for $i\geq1$, so
$(\gamma'',\delta'')$ is also strict shifted and
$l''-s''=\sum_{i=1}^{k-2}m_i+k-2.$

Repeat the above process to generate
$(\alpha^{k-2},\beta^{k-2}),\ldots,(\alpha^{1},\beta^{1})$
respectively and let $D'=D$, it's straightforward to see that the
$k$-marked Durfee symbol \break $\eta=\left(\begin{array}{cccc}
\alpha^k,&\alpha^{k-1},&\ldots,&\alpha^1\\[2pt]
\beta^k,&\beta^{k-1},&\ldots,&\beta^1
\end{array}
\right)_D$ is counted by
$\mathcal{D}^{ss}_k(m_1,m_2,\ldots,m_k;n)$. \qed

We now illustrate the reverse map $\Phi^{-1}$ by going through an
example in details. Take $m_1=1,m_2=1,m_3=0$, and let
\[\left(\begin{array}{l}
\gamma\\
\delta
\end{array}\right)_D=\left(\begin{array}{ccccccccccc}
6&6&3&3&3&3&2&2&1&1&1\\
5&5&4&2&1&1&1
\end{array}
\right)_6,\]
 we first split $(\gamma,
\delta)$ to get $(\alpha^3, \beta^3)$, note that the divisional
part $\delta_j$ is the smallest part satisfying $\delta_j\geq
\gamma_{j+1}.$
 \[
 \left(\begin{array}{l}
\alpha^3\\
\beta^3
\end{array}\right)=\left(\begin{array}{cccc}
6&6&3\\
5&5&4
\end{array}
\right),\left(\begin{array}{l}
\gamma'\\
\delta'
\end{array}\right)=\left(\begin{array}{cccccccc}
3&3&3&2&2&1&1&1\\
&2&1&1&1&&&
\end{array}
\right),\] we then split $(\gamma', \delta')$ to generate
$(\alpha^2, \beta^2)$, the divisional part $\delta'_j$ is the
smallest part satisfying $\delta'_j\geq \gamma'_{j+3}$. The
remaining part of $(\gamma', \delta')$ is just $(\alpha^1,
\beta^1)$.

\[
 \left(\begin{array}{l}
\alpha^2\\
\beta^2
\end{array}\right)=\left(\begin{array}{cccccc}
3&3&3&2&2&1\\
2&1&1&1
\end{array}
\right), \left(\begin{array}{l}
\alpha^1\\
\beta^1
\end{array}\right)=\left(\begin{array}{cc}
1&1\\
&
\end{array}
\right).
\]
Thus we get
\begin{equation}\label{example}
\eta=\left(\begin{array}{ccccccccccccc}
6_3&6_3&3_3&3_2&3_2&3_2&2_2&2_2&1_2&1_1&1_1\\
5_3&5_3&4_3&&2_2&1_2&1_2&1_2&&
\end{array}
\right)_6.\end{equation}

By Theorems \ref{r-kms-ds} and \ref{dy-op}, we then deduce the
following partition identity:
\begin{core} For $m_i\geq 0$ and $k\geq 1$,  we have
\begin{equation}
\mathcal{D}^{ss}_k(m_1,m_2,\ldots,m_k;n)=N\left(\sum_{i=1}^km_i+k-1;n\right).
\end{equation}
\end{core}

To establish the relation between $k$-marked strict shifted Durfee
symbols and $k$-marked Durfee symbols, we need define a statistic
on $k$-marked Durfee symbols. In the same way, we first define
this statistic on a pair of partitions.

 For $(\gamma,\delta)$, the part $\delta_i$ is said to be balanced if
$\gamma_{i+1}\leq \delta_i$ and the number of parts greater than
$\delta_{i}$ in $\gamma$ ($\gamma_1$ is not counted) {\it is equal
to} the number of unbalanced parts before $\delta_{i}$ in
$\delta$. For example, let
 $\left(\begin{array}{c}\gamma\\ \delta
 \end{array}
 \right)=\left(\begin{array}{ccccc}
4&3&3&1&1\\
&\underline{3}&2&2&
\end{array}\right)$, the first part  $3$ of $\delta$ is balanced while the third part $2$
is not balanced, although it satisfies the first condition
$\gamma_{i+1}\leq \delta_i$, there are two parts greater than 2
($4$ is not counted) in $\gamma$, while there is one unbalanced
part (the second part 2) before the third part 2 in $\delta$.

It should be pointed out that  for any unbalanced part $\delta_i$,
the number of parts greater than $\delta_{i}$ in $\gamma$
($\gamma_1$ is not counted) {\it is greater than} the number of
unbalanced parts before $\delta_{i}$ in $\delta$. We will state
this in the following proposition:
\begin{pro}\label{pro-bp}
Let $(\gamma,\delta)$ be a pair of partitions, for any part
$\delta_j$ of $\delta$, let $d_j$ denote the difference between
the number of parts greater than $\delta_j$ in $\gamma$
{\rm(}except for $\gamma_1${\rm)} and the number of unbalanced
parts before $\delta_j$ in $\delta$, we then have $d_j\geq 0$.
\end{pro}
For the above example, we see that $d_1=0, d_2=2, d_3=1.$

 \pf If $\gamma_{j+1}> \delta_{j}$, it's clear to see that $d_j\geq 1$; If
 $\gamma_{j+1}\leq\delta_{j}$, we consider the following two
 cases:
\begin{kst}
\item[Case 1] When $\delta_{j-1}<\gamma_{j}$, and at this time,
there are $j-1$ parts greater than $\delta_j$ in $\gamma$
($\gamma_1$ is not counted) and the number of unbalanced parts
before $\delta_j$ in $\delta$ is less than or equal to $j-1$, so
$d_j\geq 0$.

\item[Case 2] When $\delta_{j-1}\geq \gamma_{j}$, let $t$ be
largest nonnegative  integer less than $j$ such that
$\delta_{t}<\gamma_{t+1}$. Here we assume that $\delta_0=0$ and
then $t$ must exist.  From Case 1, we know that $d_{t+1}\geq 0$.
We use the induction to prove that $d_i\geq 0$ for $t+2\leq i\leq
j$.

We first prove that $d_{t+2}\geq 0$. Note that
$\gamma_{t+1}>\delta_{t+1}\geq \gamma_{t+2}$.  If
$\delta_{t+2}<\gamma_{t+2}$, then $d_{t+2}\geq d_{t+1}\geq 0$; If
$\delta_{t+2}\geq \gamma_{t+2}$, then $d_{t+2}\geq d_{t+1}-1$. In
particular, when $d_{t+1}=0$, then $\delta_{t+1}$ is balanced,
thus $d_{t+2}=d_{t+1}=0$. So $d_{t+2}\geq 0.$

Assume that $d_i\geq 0$ and $\gamma_p> \delta_i\geq \gamma_{p+1}$.
We will prove that $d_{i+1}\geq 0$. If $\delta_{i+1}\geq
\gamma_{p+1}$, the similar argument on the case for $\delta_{t+2}$
could show  that  $d_{i+1}\geq 0.$ If $\delta_{i+1}<\gamma_{p+1}$,
and then $d_{i+1}\geq d_i\geq 0$. \qed
\end{kst}

We use $b(\gamma,\delta)$ to denote the number of balanced parts
in $(\gamma,\delta).$ Clearly, $0\leq b(\gamma,\delta)\leq
l(\delta)$. Let  $\eta=\left(\begin{array}{cccc}
\alpha^k,&\alpha^{k-1},&\ldots,&\alpha^1\\[2pt]
\beta^k,&\beta^{k-1},&\ldots,&\beta^1
\end{array}
\right)_D$ be  a $k$-marked Durfee symbol, we define ${\rm
nb}_i(\eta)$, called the $i$th balanced number by
\[{\rm
nb}_i(\eta)=\left\{
\begin{array}{ll}
b(\alpha^i,\beta^i) \ \ &\text{ for }\ 1\leq i<k,\\[5pt]
0 \ \ &\text{ for }\ i=k.
\end{array}\right.
\]
For $\left(\begin{array}{cccccccc}
4_3&4_3&3_2&2_2&2_2&2_1&2_1&1_1\\
5_3&&&\underline{3_2}&\underline{2_2}&&\underline{2_1}&1_1
\end{array}\right)_5$, we have ${\rm
nb}_1=1,\,{\rm nb}_2=2,\,{\rm nb}_3=0.$

We next state a theorem concerning  strict shifted partitions.
\begin{thm}\label{lem-ssp}
Given two nonnegative integers $r,m$, there is a bijection $\psi$
between the set of pairs of partitions $(\alpha,\beta)$ of $n$
with $\beta_1\leq \alpha_1$ where there are  $r$ balanced parts
and the difference of the lengths of $\alpha$ and $\beta$ equals
to $m$ and the set of strict shifted partitions
$(\bar{\alpha},\bar{\beta})$ of $n$ where the difference of the
lengths of $\bar{\alpha}$ and $\bar{\beta}$ equals to $m+2r$.
\end{thm}
\pf {\bf The map $\psi$:} Let $(\alpha,\beta)$ be a pair of
partitions with $r$ balanced parts and $l(\alpha)-l(\beta)=m$. The
strict shifted partitions $(\bar{\alpha},\bar{\beta})$ is
constructed as follows. $\bar{\alpha}$ is composed of all parts of
$\alpha$ and all balanced parts of $\beta$. $\bar{\beta}$ consists
of all unbalanced parts of $\beta$. Take an example, let
$\left(\begin{array}{c}\alpha\\ \beta
 \end{array}
 \right)=\left(\begin{array}{cccccc}
6&5&5&3&3&2\\
&\underline{5}&4&4&\underline{3}
\end{array}\right)$, where the underlined parts in $\beta$ are
balanced. According to the above construction, we then get $\left(\begin{array}{c}\bar{\alpha}\\
\bar{\beta}
 \end{array}
 \right)=\left(\begin{array}{cccccccc}
6&\underline{5}&5&5&\underline{3}&3&3&2\\
&4&4
\end{array}\right)$. It's clear to see that
$|\alpha|+|\beta|=|\bar{\alpha}|+|\bar{\beta}|$,
$l(\bar{\alpha})-l(\bar{\beta})=l(\alpha)+r-(l(\beta)-r)=m+2r.$ By
Proposition \ref{pro-bp}, one can also easily know that
$(\bar{\alpha},\bar{\beta})$ is strict shifted.

{\bf The reverse map $\psi^{-1}$:} Let
$(\bar{\alpha},\bar{\beta})$ be a strict shifted partition where
the difference  of the lengths of $\bar{\alpha}$ and $\bar{\beta}$
is $m+2r$, that is $l(\bar{\alpha})-l(\bar{\beta})=m+2r$. We now
construct a pair of partitions $(\alpha,\beta)$ with $r$ balanced
parts and $l(\alpha)-l(\beta)=m$.

First of all, attach subscript $g_i$ for each part
$\bar{\alpha}_i$ of $\bar{\alpha}$, where $g_i$ denotes the
difference between the number of parts before $\bar{\alpha}_i$ in
$\bar{\alpha}$ ($\bar{\alpha}_1$ is not counted) and the number of
parts greater than or equal to $\bar{\alpha}_i$ in $\beta$. We let
$g_1=0$.

For example,
if $\left(\begin{array}{c}\bar{\alpha}\\
\bar{\beta}
 \end{array}
 \right)=\left(\begin{array}{ccccccccc}
6&5&5&5&3&3&3&2&1\\
&4&4&3
\end{array}\right)$, attach the subscripts for all parts of $\bar{\alpha}$ to
get
$\left(\begin{array}{ccccccccc}
6_0&5_0&5_1&5_2&3_0&3_1&3_2&2_3&1_4\\
&4&4&3
\end{array}\right).$

One could easily know that $g_2=0$ and $g_i\geq 0$ for any $3\leq
i\leq l(\bar{\alpha})$ from the fact that
$(\bar{\alpha},\bar{\beta})$ is strict shifted. Let
$\bar{\alpha}_{t_i}$ be the smallest part in all of parts of
$\alpha$ with subscript equal to $i$, we have the following
conclusion:

\begin{lem}\label{temp}
For $0\leq i \leq m+2r-2$, $\bar{\alpha}_{t_i}$ exists, and
$\bar{\alpha}_{t_0}\geq \bar{\alpha}_{t_1}\geq \cdots \geq
\bar{\alpha}_{t_{m+2r-2}}$.
\end{lem}
In the above example, $m+2r=6$ and
$\bar{\alpha}_{t_0}=3,\bar{\alpha}_{t_1}=3,\bar{\alpha}_{t_2}=3,\bar{\alpha}_{t_3}=2,\bar{\alpha}_{t_4}=1$.

  \pf We use the induction to show that the
sequence of subscripts $\{g_1,g_2,\ldots,g_{l(\bar{\alpha})}\}$
consists of all nonnegative integers less than $m+2r-1$.
Obviously, $0$ is in this sequence. Assume that $i$ is in this
sequence, that is there is a part $\bar{\alpha}_p$ such that
$g_p=i$, we now prove that $i+1$ is also in this sequence. By the
induction hypothesis, we know the subscript of the part
$\bar{\alpha}_p$ is $i$, that is
$\bar{\beta}_{p-i-1}<\bar{\alpha}_p\leq \bar{\beta}_{p-i-2}$. Let
$l(\bar{\beta})=s$ and note that $i\leq m+2r-3$, we have
$s+i+3\leq s+m+2r=l(\bar{\alpha})$. If $\bar{\alpha}_{s+i+3}\leq
\bar{\beta}_{s}$, then the subscript of $\bar{\alpha}_{s+i+3}$ is
$i+1$; Otherwise, there must exist $p+1\leq j\leq s+i+2$ such that
$\bar{\beta}_{j-i-2}<\bar{\alpha}_{j}\leq \bar{\beta}_{j-i-3}$
this is because that $\bar{\alpha}_{p+1}\leq \bar{\beta}_{p-i-2}$
and $\bar{\alpha}_{s+i+2}>\bar{\beta}_s$. Hence, the subscript of
$\bar{\alpha}_{j}$ is $i+1$. Therefore, $\bar{\alpha}_{t_i}$
exists for $0\leq i \leq m+2r-2$ and $\bar{\alpha}_{t_0}\geq
\bar{\alpha}_{t_1}\geq \cdots \geq \bar{\alpha}_{t_{m+2r-2}}$ when
note that given a part $\bar{\alpha}_p$ with subscript $i$, we
could always find a part after $\bar{\alpha}_p$ whose subscript is
$i+1$. \qed

Let $\gamma$ be a partition having $r$ parts whose parts are
 $\bar{\alpha}_{t_0},
\bar{\alpha}_{t_2},\ldots,\bar{\alpha}_{t_{r-1}}$ respectively.
Take $r=2$ in the above example, $\gamma=(3,3)$. We now construct
the partitions $(\alpha,\beta)$. The partition $\alpha$ consists
of all parts in $\bar{\alpha}$, while not in $\gamma$. $\beta$ is
composed of all parts both in $\bar{\beta}$ and $\gamma$.
In the above example, we therefore get $\left(\begin{array}{c}\alpha\\
\beta
 \end{array}
 \right)=\left(\begin{array}{ccccccccc}
6&5&5&5&3&2&1\\
&4&4&3&\underline{3}&\underline{3}
\end{array}\right)$.

It's obvious to see that
$|\bar{\alpha}|+|\bar{\beta}|=|\alpha|+|\beta|$,
$l(\alpha)-l(\beta)=l(\bar{\alpha})-l(\gamma)-[l(\bar{\beta})+l(\gamma)]=m$.
We now show that $(\alpha,\beta)$ has exactly $r$ balanced parts.
From Lemma \ref{temp} and the definition of $\gamma$, we know that
for each  part $\beta_t$ from $\gamma$ in $\beta$, the number of
parts greater than $\beta_t$ in $\alpha$ equals the number of
parts from $\bar{\beta}$ greater than or equal to $\beta_t$ in
$\beta$. Thus we just need to prove that the parts from
$\bar{\beta}$ in $\beta$ are unbalanced. We use induction on the
part from $\bar{\beta}$ in $\beta$. We first verify the largest
part $\bar{\beta}_1$ of $\bar{\beta}$ is unbalanced. Supposed that
there are $t$ parts $\gamma_1,\gamma_2,\ldots,\gamma_t$ from
$\gamma$ greater than $\bar{\beta}_1$, then
$\gamma_i=\bar{\alpha}_{i+1},i=1,2,\ldots,t$ and
$\alpha_2=\bar{\alpha}_{t+2}$. We claim that
$\alpha_2=\bar{\alpha}_{t+2}>\bar{\beta}_1$. Recall that the part
$\gamma_i$ is the smallest part whose subscript is $i$. If
$\alpha_2=\bar{\alpha}_{t+2}\leq \bar{\beta}_1$, the subscript of
$\bar{\alpha}_{t+2}$ is less than $t+1$, this contradicts to the
definition of $\gamma_i$. Clearly, these $t$ parts from $\gamma$
are balanced, and $\bar{\beta}_1$ is not balanced. We now consider
the part $\bar{\beta}_j$ from $\bar{\beta}$ in $\beta$. Assume
that all parts from $\bar{\beta}$ before $\bar{\beta}_j$ in
$\beta$ are not balanced and there are $t$ parts from $\gamma$
before $\bar{\beta}_j$. We next justify
$\alpha_{j+1}>\bar{\beta}_j$ and then by the hypothesis, we know
that there are $j-1$ unbalanced parts before $\bar{\beta}_j$,
while there are at least $j$ parts larger than $\bar{\beta}_j$ in
$\alpha$, so $\bar{\beta}_j$ is unbalanced. Since there are $t$
parts from $\gamma$ before $\bar{\beta}_j$, then
$\alpha_{j+1}=\bar{\alpha}_{t+j+1}$ and if
$\bar{\alpha}_{t+j+1}=\alpha_{j+1}\leq \bar{\beta}_j$, then the
subscript of $\bar{\alpha}_{t+j+1}$ is less than $t$, which
contradicts to the definition of the parts of $\gamma$, so
$\alpha_{j+1}>\bar{\beta}_j$ and we therefore complete the proof.
\qed

The next theorem gives  a relation between $k$-marked strict
shifted Durfee symbols and $k$-marked Durfee symbols.

\begin{thm}\label{r-km-kms}
Given $2k$ nonnegative integers $m_1,m_2,\ldots,m_k,$ and
$t_1,t_2,\ldots,t_k$ where $t_k=0$. There is a bijection $\Psi$
between the set of $k$-marked Durfee symbols of $n$ with $i$th
rank equal to $m_i$ and $i$th balanced number equal to $t_i$ and
the set of $k$-marked strict shifted Durfee symbols of $n$ with
$i$th rank equal to $m_i+2t_i$.
\end{thm}
\pf Let $\eta=\left(\begin{array}{cccc}
\alpha^k,&\alpha^{k-1},&\ldots,&\alpha^1\\[2pt]
\beta^k,&\beta^{k-1},&\ldots,&\beta^1
\end{array}
\right)_D$ be a $k$-marked Durfee symbol with $i$th rank equal to
$m_i$ and $i$th balanced number
 equal to $t_i$. We now apply the
bijection $\psi$ in Theorem \ref{lem-ssp} on each vector
$(\alpha^i,\beta^i)$ of $\eta$ except for $k$th vector
$(\alpha^k,\beta^k)$, to generate
$(\bar{\alpha}^i,\bar{\beta}^i)$. From Theorem \ref{lem-ssp}, we
know that $(\bar{\alpha}^i,\bar{\beta}^i)$ is strict shifted and
$l(\bar{\alpha}^i)-l(\bar{\beta}^i)=l(\alpha^i)-l(\beta^i)+2t_i$.
Let $\bar{\eta}=\left(\begin{array}{cccc}
\alpha^k,&\bar{\alpha}^{k-1},&\ldots,&\bar{\alpha}^1\\[2pt]
\beta^k,&\bar{\beta}^{k-1},&\ldots,&\bar{\beta}^1
\end{array}
\right)_D$ which has the same subscript and the same $k$th vector
with $\eta$. It's obvious to see that $\bar{\eta}$ is a $k$-marked
strict shifted Durfee symbol with $i$th rank equal to $m_i+2t_i$.
\qed

By Theorem \ref{r-km-kms}, one can derive the following identity
readily.
\begin{core}\label{c-km-kms}
For $m_i\geq 0$ and $k\geq 2$, we have
\begin{equation}
\mathcal{D}_k(m_1,m_2,\ldots,m_k;n)=\sum_{t_1,\ldots,\,t_{k-1}=0}^{{+\infty}}\mathcal{D}^{ss}_k(m_1+2t_1,\ldots,m_{k-1}+2t_{k-1},m_k;n).
\end{equation}
\end{core}
Combine  Corollaries \ref{r-kms-ds} and \ref{c-km-kms} to get:
\begin{thm}\label{c-km-ds-1}
For $m_i\geq 0$ and $k\geq 2$, we have
\begin{equation}
\mathcal{D}_k(m_1,m_2,\ldots,m_k;n)=\sum_{t_1,\ldots,\,t_{k-1}=0}^{+\infty}
N\left(\sum_{i=1}^km_i+2\sum_{i=1}^{k-1}t_i+k-1;n\right).
\end{equation}
\end{thm}

The following compact form  of Theorem \ref{c-km-ds-1} can be
easily obtained upon utilizing the fact that the number of
solutions to $t_1+t_2+\cdots+t_{k-1}=j$ in nonnegative integers is
$j+k-2 \choose k-2$.
\begin{thm}\label{c-km-ds-2}
For $m_i\geq 0$ and $k\geq 2$, we have
\begin{equation}
\mathcal{D}_k(m_1,m_2,\ldots,m_k;n)=\sum_{j=0}^{+\infty}{j+k-2
\choose k-2}N\left(\sum_{i=1}^km_i+2j+k-1;n\right).
\end{equation}
\end{thm}

We next generalize  Theorem \ref{c-km-ds-2} to give Theorem
\ref{main} which holds for any integer $m_i$.  To do this, we
prove the following conclusion by constructing a simple bijection
$\Theta$.
\begin{thm} \label{r-sym}
For $k\geq 1$ and $1\leq p\leq k$, we have
\begin{equation}
\mathcal{D}_k(m_1,\ldots,m_p,\ldots
m_k;n)=\mathcal{D}_k(m_1,\ldots,-m_p,\ldots m_k;n).
\end{equation}
\end{thm}
\pf  Let $\eta=\left(\begin{array}{cccc}
\alpha^k,&\alpha^{k-1},&\ldots,&\alpha^1\\[2pt]
\beta^k,&\beta^{k-1},&\ldots,&\beta^1
\end{array}
\right)_D$ be a $k$-marked Durfee symbol with $i$th rank equal to
$m_i$. We will construct another $k$-marked Durfee symbol
$\bar{\eta}$ with $i$th rank equal to $\bar{m}_i$ such that
$\bar{m}_p=-m_p$ and $\bar{m}_i=m_i$ for $i\neq p$.

Define \[\bar{\eta}=\left(\begin{array}{ccccccc}
\alpha^k,&\ldots,&\alpha^{p+1}&\bar{\alpha}^p,&\alpha^{p-1},&\ldots,&\alpha^1\\[2pt]
\beta^k,&\ldots,&\beta^{p+1},&\bar{\beta}^p,&\beta^{p-1},&\ldots,&\beta^1
\end{array}
\right)_D,\] where $\bar{\alpha}^k=\beta^k$ and
$\bar{\beta}^k=\alpha^k$ for $p=k$; When $p\neq k$,
$\bar{\alpha}^p$ consists of all parts of $\beta^p$ and the
largest part $\alpha^p_1$ of $\alpha^p$. $\bar{\beta}^p$ consists
of all parts of $\alpha^p$ except for the largest part
$\alpha^p_1$. It's clear to see that
$\bar{m}_p=l(\bar{\alpha}^p)-l(\bar{\beta}^p)-1=l(\beta^p)+1-[l(\alpha^p)-1]-1=-[l(\alpha^p)-l(\beta^p)-1)]=-m_p$
for $p\neq k$ and
$\bar{m}_k=l(\bar{\alpha}^k)-l(\bar{\beta}^k)=-[l(\alpha^k)-l(\beta^k)]=-m_k,$
so $\bar{\eta}$ is desired. \qed

By Theorem \ref{r-sym}, we could generalize Theorem
\ref{c-km-ds-1} to give the following theorem which is useful to
prove a relationship between $k$-marked Durfee symbols and the
symmetrized rank moment given by Andrews (see Theorem \ref{int}).

\begin{thm}\label{c-km-ds-1'}
For $k\geq 2$, we have
\begin{equation}
\mathcal{D}_k(m_1,m_2,\ldots,m_k;n)=\sum_{t_1,\ldots,\,t_{k-1}=0}^{+\infty}
N\left(\sum_{i=1}^k|m_i|+2\sum_{i=1}^{k-1}t_i+k-1;n\right).
\end{equation}
\end{thm}
Theorem \ref{main} is a compact form of Theorem \ref{c-km-ds-1'}
which immediately follows from Theorems \ref{c-km-ds-2} and
\ref{r-sym}.

\section{Andrews' identities on $k$-marked Durfee symbols }

In this section, we aim to  show  the identities  on $k$-marked
Durfee symbols given by Andrews with the help of Theorem
\ref{main}. Recall that $\mathcal{D}_k(m_1,m_2,\ldots,m_k;n)$
denotes the number of $k$-marked Durfee symbols of $n$ with $i$th
rank equal to $m_i$. Andrews considered the following generating
function for $\mathcal{D}_k(m_1,m_2,\ldots,m_k;n)$:
\begin{align}\label{def-g-km}
R_k(x_1,\ldots,x_k;q)=\sum_{m_1,\ldots,m_k=-\infty}^{+\infty}\sum_{n=0}^{+\infty}
\mathcal{D}_k(m_1,\ldots,m_k;n)x_1^{m_1}\cdots x_k^{m_k}q^n,\,
|q|<1.
\end{align}
By  applying the $k$-fold generalization of Watson's
transformation between a very-well-poised $_8\phi_7$-series and a
balanced $_4\phi_3$-series \cite[p.199, Theorem 4]{Andrews-1975},
Andrews gave the generating function $R_k(x_1,x_2,\ldots,x_k;q)$
in the following theorem.

\begin{thm}[Corollary 11, Andrews \cite{Andrews-07-a}]\label{gf-km}
\begin{align}\label{gf-km-e}
&R_k(x_1,x_2,\ldots,x_k;q)
=\frac{1}{(q;q)_\infty}\sum_{n=1}^{+\infty}(-1)^{n-1}q^{3n(n-1)/2+kn}
\frac{(1+q^n)(1-q^n)^2}{\prod_{j=1}^k(1-x_jq^n)(1-\frac{q^n}{x_j})}.
\end{align}
\end{thm}
\pf We will reformulate this identity  as the partition identity
\eqref{main-e} in Theorem \ref{main}. The key step is to give a
partition interpretation of the right side hand of
\eqref{gf-km-e}. We will show that it is the generating function
for the summation on the right side of \eqref{main-e}.

First, the right hand side of \eqref{gf-km-e} can be written as
the difference of the following two terms:
\begin{align*}
&\frac{1}{(q;q)_\infty}\sum_{n=1}^{+\infty}(-1)^{n}q^{n(3n-1)/2+kn}
\frac{(1-x_1)(1-x_1^{-1})(1+q^n)}{(1-x_1q^n)(1-x_1^{-1}q^n)}\times\frac{1}{\prod_{j=2}^k(1-x_jq^n)(1-\frac{q^n}{x_j})}\nonumber \\
&\hskip
-0.5cm-\frac{1}{(q;q)_\infty}\sum_{n=1}^{+\infty}(-1)^{n}q^{n(3n-1)/2+(k-1)n}
\frac{1+q^n}{\prod_{j=2}^k(1-x_jq^n)(1-\frac{q^n}{x_j})}.\nonumber
\end{align*}
We next expand each term of the above two terms, note that
\begin{align*}
\frac{(1-x_1)(1-x_1^{-1})(1+q^n)}{(1-x_1q^n)(1-x_1^{-1}q^n)}&=\frac{(1-x_1)}{(1-x_1q^n)}
+\frac{(1-x_1^{-1})}{(1-x_1^{-1}q^n)}
=2+\sum_{\stackrel{m_1=-\infty}{m_1\neq
0}}^{{+\infty}}x_1^{m_1}[q^{n|m_1|}-q^{n(|m_1|-1)}],\\
\frac{1}{(1-x_jq^n)(1-x_j^{-1}q^n)}&=\sum_{a=0}^{+\infty}
x_j^aq^{na}\sum_{b=0}^{+\infty}
x_j^{-b}q^{nb}=\sum_{m_j=-\infty}^{{+\infty}}x_j^{m_j}\sum_{t_j=0}^{{+\infty}}q^{n(|m_j|+2t_j)}.
\end{align*}
Given $k$ integers $m_1\ldots,m_k$, it's clear to see that the
coefficients of $x_1^{m_1}x_2^{m_2}\cdots x_k^{m_k}$ on the series
expansion of  the right  hand side of \eqref{gf-km-e} are
\begin{align*}
&\frac{1}{(q;q)_\infty}\sum_{n=1}^{+\infty}(-1)^{n}q^{n(3n-1)/2+kn}[q^{n|m_1|}-q^{n(|m_1|-1)}]
\sum_{t_2,\ldots,t_k=0}^{{+\infty}}q^{n(\sum_{i=2}^k|m_i|+2\sum_{i=2}^kt_i)}\\
&=\frac{1}{(q;q)_\infty}\sum_{n=1}^{+\infty}(-1)^{n-1}q^{n(3n-1)/2+kn}(1-q^n)q^{n(|m_1|-1)}
\sum_{j=0}^{{+\infty}}{j+k-2\choose
k-2}q^{n(\sum_{i=2}^k|m_i|+2j)}\\
&=\sum_{j=0}^{{+\infty}}{j+k-2\choose
k-2}\frac{1}{(q;q)_\infty}\sum_{n=1}^{+\infty}(-1)^{n-1}q^{n(3n-1)/2}(1-q^{n})q^{n(\sum_{i=1}^k|m_i|+2j+k-1)}\\
&=\sum_{j=0}^{{+\infty}}{j+k-2\choose k-2}\sum_{n=0}^{+\infty}
N(\sum_{i=1}^k|m_i|+2j+k-1;n)q^n\\
&=\sum_{n=0}^{+\infty} q^n[\sum_{j=0}^{{+\infty}}{j+k-2\choose
k-2}N(\sum_{i=1}^k|m_i|+2j+k-1;n)],
\end{align*}
where the penultimate identity  follows from Theorem \ref{gf-r}
and we then obtain the following combinatorial interpretation:
\begin{align}\label{com-int}
&\frac{1}{(q;q)_\infty}\sum_{n=1}^{+\infty}(-1)^{n-1}q^{3n(n-1)/2+kn}
\frac{(1+q^n)(1-q^n)^2}{\prod_{j=1}^k(1-x_jq^n)(1-\frac{q^n}{x_j})}\nonumber \\[2pt]
&\hskip
-0.5cm=\sum_{m_1,\ldots,m_k=-\infty}^{+\infty}\sum_{n=0}^{+\infty}
\left[\sum_{j=0}^{+\infty}{j+k-2 \choose
k-2}N\left(\sum_{i=1}^k|m_i|+2j+k-1;n\right)\right]x_1^{m_1}\cdots
x_k^{m_k}q^n.
\end{align}
Combining  \eqref{def-g-km} and \eqref{com-int},  we reach our
conclusion that the identity \eqref{gf-km-e} can be restated as
the partition identity \eqref{main-e}. Thus, we have obtained a
combinatorial proof of \eqref{gf-km-e} based on Theorem
\ref{main}. \qed

Recently, Bringmann, Lovejoy, and Osbur defined a two-parameter
generalization of $k$-marked Durfee symbols in
\cite{Bringmann-lovejoy-osbur}. They deduced  the generating
function \cite[Theorem 2.2]{Bringmann-lovejoy-osbur} for the
two-parameter generalization of $k$-marked Durfee symbols using
the similar argument of Andrews, which reduces to the identity
\eqref{gf-km-e} when $d=e=0.$

From the generating function $R_k(x_1,x_2,\ldots,x_k;q)$ in
Theorem \ref{gf-km}, Andrews immediately found the following
symmetry of $k$-marked Durfee symbols.
\begin{thm}[Corollary 12, Andrews \cite{Andrews-07-a}]\label{sym-mds}
$\mathcal{D}_k(m_1,m_2,\ldots,m_k;n)$ is symmetric in
$m_1,\,m_2,\ldots,m_k$.
\end{thm}
\pf This symmetry can also immediately  follow from Theorem
\ref{main}.
 \qed

In fact, the composite of the bijections on Section 2 provides a
bijection for this symmetry. We take an example to explain this
process. Let
\begin{equation*}
\eta=\left(\begin{array}{ccccccccccccc}
6_3&3_2&3_2&2_2&2_2&1_2&1_1\\
5_3&3_2&3_2&1_2&&&1_1&1_1
\end{array}
\right)_6 \text{counted by }
\mathcal{D}_3(-2,1,0;68),\end{equation*} we aim to construct a
$3$-marked Durfee symbol $\bar{\eta}$ counted by
$\mathcal{D}_3(1,-2,0;68)$. We will first combine all subscripts
of $k$-marked Durfee symbol $\eta$ to get a Durfee symbol, and
then split this Durfee symbol over again to get our desired
$k$-marked Durfee symbol $\bar{\eta}$.

First, applying the bijection $\Theta$ in Theorem \ref{r-sym} into
$\eta$ to get $\eta'$ enumerated by $\mathcal{D}_3(2,1,0;68)$,
\begin{equation*}
\eta'=\left(\begin{array}{ccccccccccccc}
6_3&3_2&3_2&2_2&2_2&1_2&1_1&1_1&1_1\\
5_3&&\underline{3_2}&\underline{3_2}&1_2&
\end{array}
\right)_6, \end{equation*} we now  utilize the bijection $\Psi$ in
Theorem \ref{r-km-kms} on $\eta'$ to get a $k$-marked strict
shifted Durfee symbol. Observe that there are two balanced parts
in the second vector of $\eta'$, and there is no balanced part in
other vectors of $\eta'$. So we will get a $k$-marked strict
shifted Durfee symbol $\eta''$ which counted by
$\mathcal{D}^{ss}_3(2,1+2\times2,0;68)$,
\begin{equation*}
\eta''=\left(\begin{array}{ccccccccccccc}
6_3&3_2&3_2&3_2&3_2&2_2&2_2&1_2&1_1&1_1&1_1\\
5_3&&   1_2
\end{array}
\right)_6, \end{equation*} applying the bijection $\Phi$ in
Theorem \ref{r-kms-ds} to get $\eta'''$ which counted by
$\mathcal{D}_1(9;68)$
\begin{equation*}
\eta'''=\left(\begin{array}{ccccccccccccc}
6&3&3&3&3&2&2&1&1&1&1\\
5&1
\end{array}
\right)_6. \end{equation*} Thus we complete the first step. We
next split the Durfee symbol $\eta'''$ over again.

First, apply the reverse map $\Phi^{-1}$ in Theorem \ref{r-kms-ds}
on $\eta'''$ to get $\bar{\eta}''$ which counted by
$\mathcal{D}^{ss}_3(1,2+2\times2,0;68)$,
\begin{equation*}
\bar{\eta}''=\left(\begin{array}{ccccccccccccc}
6_3&3_2&\underline{3_2}&\underline{3_2}&3_2&2_2&2_2&1_2&1_2&1_1&1_1\\
5_3&1_2
\end{array}
\right)_6, \end{equation*} using the reverse map $\Psi^{-1}$ in
Theorem \ref{r-km-kms} on $\bar{\eta}''$, we get $k$-marked Durfee
symbol $\bar{\eta}'$ counted by $\mathcal{D}_3(1,2,0;68)$
\begin{equation*}
\bar{\eta}'=\left(\begin{array}{ccccccccccccc}
6_3&3_2&3_2&2_2&2_2&1_2&1_2&1_1&1_1\\
5_3&&3_2&3_2&1_2
\end{array}
\right)_6, \end{equation*} Finally, we obtain the desired
$k$-marked Durfee symbol $\bar{\eta}$ counted by
$\mathcal{D}_3(1,-2,0;68)$ when applying the bijection $\Theta$ in
Theorem \ref{r-sym} on $\bar{\eta}'$.
\begin{equation*}
\bar{\eta}=\left(\begin{array}{ccccccccccccc}
6_3&3_2&3_2&3_2&1_2&&1_1&1_1\\
5_3&3_2&2_2&2_2&1_2&1_2
\end{array}
\right)_6.\end{equation*}

By the generating function $R_k(x_1,x_2,\ldots,x_k;q)$ and the
generating function for $\eta_{2k}(n)$, Andrews  showed that the
number of $(k+1)$-marked Durfee symbols of $n$ equals the
symmetrized $(2k)$-th moment function at $n$ in
\cite{Andrews-07-a}, that is

\begin{thm}[Corollary 13, Andrews \cite{Andrews-07-a}]\label{int}
For $k\geq 1$,
\begin{equation}\label{int-e}
\mathcal{D}_{k+1}(n)=\eta_{2k}(n).
\end{equation}

\end{thm}

\pf Recall that
\begin{align*}
\eta_{2k}(n)&=\sum_{m=-\infty}^{{+\infty}}{m+k-1 \choose
2k}N(m;n)=\sum_{m=1}^{{+\infty}}\left[{m+k-1 \choose 2k}+{m+k
\choose 2k}\right]N(m;n),
\end{align*}
where the second equality follows from the
 rank symmetry $N(-m;n)=N(m;n)$ and the fact ${-m+k-1 \choose 2k}={m+k \choose 2k}$.
\begin{align*}
\mathcal{D}_{k+1}(n)&=\sum_{m_1,\ldots,m_{k+1}=-\infty}^{+\infty}
\mathcal{D}_{k+1}(m_1,m_2,\ldots,m_{k+1};n)\\
&=\sum_{m_1,\ldots,m_{k+1}=-\infty}^{+\infty}\sum_{t_1,\ldots,\,t_{k}=0}^{+\infty}
N\left(\sum_{i=1}^{k+1}|m_i|+2\sum_{i=1}^{k}t_i+k;n\right),\end{align*}
where the second equality follows from  Theorem \ref{c-km-ds-1'}.
So it suffices to show that the number of solutions to
$|m_1|+\cdots+|m_{k+1}|+2t_1++\cdots+2t_k=m-k$ where $m_i$ is
integer, and $t_i$  is nonnegative integer equals to ${m+k-1
\choose 2k}+{m+k \choose 2k}$.

Let $c(n)$ denote the number of solutions to
$|m_1|+|m_2|+\cdots+|m_{k+1}|+2t_1+2t_2+\cdots+2t_k=n$ where $m_i$
is integer and $t_i$ is nonnegative integer. It's easy to know
that the generating function for $c(n)$ is
\begin{align*}
\sum_{n=0}^{{+\infty}}c(n)q^n&=(1+2q+2q^2+2q^3+\cdots)^{k+1}(1+q^2+q^4+q^6+\cdots)^k \\
&=\frac{(1+q)^{k+1}}{(1-q)^{k+1}}\times\frac{1}{(1-q^2)^k} \\
&=(1+q)\times\frac{1}{(1-q)^{2k+1}}\\
&=(1+q)\sum_{n=0}^{{+\infty}}{2k+n \choose n}q^n\\
&=\sum_{n=0}^{{+\infty}}{2k+n \choose
2k}q^n+\sum_{n=1}^{{+\infty}}{2k+n-1 \choose 2k}q^n.
\end{align*}
Comparing coefficients of $q^n$ in the above expression, we obtain
\[c(n)={2k+n \choose
2k}+{2k+n-1 \choose 2k}.\] Thus we reach our conclusion. \qed

By partial fraction expansion, Andrews \cite{Andrews-07-a} also
gave the following relationship between the generating function
for $k$-marked Durfee symbols and the generating function for
Durfee symbols, which plays an important role in the study of
Ramanujan-type congruences for $k$-marked Durfee symbols.

\begin{thm}[Theorem 7, Andrews \cite{Andrews-07-a}]\label{thm7}
\begin{align}\label{thm-7}
R_k(x_1,x_2,\ldots,x_k;q)=\sum_{i=1}^k\frac{R_1(x_i;q)}{\prod_{\stackrel{j=1}{
j\neq i}}^k(x_i-x_j)(1-\frac{1}{x_ix_j})}.
\end{align}
\end{thm}
\pf Similarly, we will restate this identity as the partition
identity \eqref{main-e} in Theorem \ref{main}.  We first consider
the series expansion of  the right hand side of \eqref{thm-7}.  To
do this, we need to work in a larger ring: the field of iterated
Laurent series $K\ll
x_k,x_{k-1},\ldots,x_1\gg=K((x_k))((x_{k-1}))\cdots((x_1))$ where
$K=\mathbb{C}(q)$, in which all series are regarded first as
Laurent series in $x_1$, then as Laurent series in $x_2$, and so
on. For more detailed account of the properties of this field,
with other applications, see \cite{Xin-04} and \cite{Xin-05}.

Every element of $K\ll x_k,x_{k-1},\ldots,x_1\gg$ has a unique
iterated Laurent series expansion. In particular, the series
expansion of $1/(1-\frac{1}{x_ix_j})$ is:
\begin{equation}
\frac{1}{1-\frac{1}{x_ix_j}}=-\frac{x_ix_j}{1-x_ix_j}=-\sum_{l=1}^{+\infty}x_i^{l}x_j^{l}.
\end{equation}
The series expansions of $1/(x_i-x_j)$ will be especially
important. If $j<i$, then
\begin{equation*}
\frac{1}{x_i-x_j}=\frac{x^{-1}_i}{1-\frac{x_j}{x_i}}=\sum_{l=0}^{+\infty}x_i^{-l-1}x_j^{l}.
\end{equation*}
However, if $j>i$ then this expansion is not valid and instead we
have the expansion:
\begin{equation*}
\frac{1}{x_i-x_j}=\frac{-x^{-1}_j}{1-\frac{x_i}{x_j}}=-\sum_{l=0}^{+\infty}x_i^{l}x_j^{-l-1}.
\end{equation*}
Thus for $j<i$, the series expansion of
$\frac{1}{(x_i-x_j)(1-x^{-1}_ix^{-1}_j)}$ is
\begin{equation}\label{serexp-1}
\frac{1}{(x_i-x_j)(1-x^{-1}_ix^{-1}_j)}=-\sum_{l=1}^{+\infty}x_i^{l}x_j^{l}\sum_{m=0}^{+\infty}x_i^{-m-1}x_j^{m}
=-\sum_{m_j=1}^{+\infty}x_j^{m_j}\sum_{t_j=0}^{m_j-1}x_i^{2t_j-m_j+1},
\end{equation}
and for $j>i$, we have the following expansion:
\begin{equation}\label{serexp-2}
\frac{1}{(x_i-x_j)(1-x^{-1}_ix^{-1}_j)}=\sum_{l=1}^{+\infty}x_i^{l}x_j^{l}\sum_{m=0}^{+\infty}x_i^{m}x_j^{-m-1}
=\sum_{m_j=-\infty}^{+\infty}x_j^{m_j}\sum_{t_j=0}^{+\infty}x_i^{|m_j|+2t_j+1}.
\end{equation}
We now consider the series expansion of the $i$th term of the
right hand side of \eqref{thm-7}.
\begin{equation}\label{ithterm-i}
\frac{R_1(x_i;q)}{\prod_{\stackrel{j=1}{ j\neq
i}}^k(x_i-x_j)(1-\frac{1}{x_ix_j})}=\frac{\sum_{n=0}^{+\infty}\sum_{m_i=-\infty}^{+\infty}\mathcal{D}_1(m_i;n)x_i^{m_i}q^n}{\prod_{\stackrel{j=1}{
j\neq i}}^k(x_i-x_j)(1-\frac{1}{x_ix_j})}.
\end{equation}
 Observe that the numerator in the above term is a series expansion in $x_i$ and
 $q$ and by the series expansions \eqref{serexp-1} and \eqref{serexp-2},  we obtain a series expansion of \eqref{ithterm-i},
 in which the exponents of $x_1,x_2,\ldots,x_{i-1}$
must be positive  and the coefficients of $x_1^{m_1}\cdots
x_k^{m_k}q^n$ for $m_1,\ldots,m_{i-1}\geq1$ are
\begin{align} \label{coef-iterm}
(-1)^{i-1}\sum_{\mathbf{t}_i}\mathcal{D}_1
\left(\sum_{j=1}^{i-1}m_j+m_i-\sum_{j=i+1}^k|m_j|-2|\mathbf{t}_i|-(k-1);n\right),
\end{align}
where the sum ranges over all sequences
$\mathbf{t}_i=(t_1,\ldots,\hat{t}_i,\ldots,t_k)$ (omitting $t_i$)
where $0\leq t_j<|m_j|$ for $j<i$ and $t_j$ could be arbitrary
nonnegative integer for $j>i$. Define
$|\mathbf{t}_i|=\sum_{\stackrel{j=1}{j\neq i}}^kt_j.$

Thus, we obtain a series expansion of the right hand of
\eqref{thm-7}:
\begin{align}\label{rhand}
&\sum_{i=1}^k\frac{R_1(x_i;q)}{\prod_{\stackrel{j=1}{ j\neq
i}}^k(x_i-x_j)(1-\frac{1}{x_ix_j})}=\sum_{i=1}^{k}
\sum_{m_1,\ldots,m_{i-1}=1}^{+\infty}
\sum_{m_i,\ldots,m_{k}=-\infty}^{+\infty}\sum_{n=0}^{+\infty}x_1^{m_1}\cdots
x_{k}^{m_k}q^n \\
 &\hskip 1cm \times \left[(-1)^{i-1}\sum_{\mathbf{t}_i}\mathcal{D}_1
\left(\sum_{j=1}^{i-1}m_j+m_i-\sum_{j=i+1}^k|m_j|-2|\mathbf{t}_i|-(k-1);n\right)\right]\nonumber.
\end{align}

Let $\mathbf{m}_i=(m_1,\ldots,m_k)$ where  $m_j\geq 1$  for $j<i$,
$m_i\leq 0$, and others could be arbitrary integers. Define
$\mathbf{x}^{\mathbf{m}_i}=x_1^{m_1}\cdots x_k^{m_k}$. Obviously,
the term $\mathbf{x}^{\mathbf{m}_i}q^n$ would be appeared in the
series expansions of the first $i$ terms of \eqref{rhand}. We next
use the induction to prove the coefficients of
$\mathbf{x}^{\mathbf{m}_i}q^n$ in the series expansion of
\eqref{rhand} equal to
\begin{equation}\label{rghthandthm7}
\sum_{j=0}^{+\infty}{j+k-2 \choose
k-2}N\left(\sum_{i=1}^k|m_i|+2j+k-1;n\right).
\end{equation}

 Let $\daleth_i(m_1,\ldots,m_{i-1})$ denote the set
of all sequences $\mathbf{t}_i$ of nonnegative integers (omitting
the $i$th vector $t_i$) such that  $t_j$  less than $|m_j|$ for
$j<i$. The following two lemmas are useful in our argument.

\begin{lem}\label{seq-lem-2}
\begin{itemize}
\item[(1)]The  number of sequences $\mathbf{t}'_{p+1} \in
\daleth_{p+1}(m_1,\ldots,m_{p}) $  equals the number of sequences
$\mathbf{t}_p \in \daleth_p(m_1,\ldots,m_{p-1})$ where the
$(p+1)$th vector $t_{p+1}<|m_p|$ such that
$|\mathbf{t}_p|=|\mathbf{t}'_{p+1}|$.

\item[(2)]The number of sequences $\mathbf{t}_p \in
\daleth_p(m_1,\ldots,m_{p-1})$ where the $(p+1)$th vector
$t_{p+1}\geq |m_p|$ is equal to the  number of sequences
$\mathbf{t}'_{p}\in \daleth_p(m_1,\ldots,m_{p-1})$ such that
$|\mathbf{t}_p|=|\mathbf{t}'_{p}|+|m_p|$.
\end{itemize}
\end{lem}

\pf Given a sequence
$\mathbf{t}_p=(t_1,\ldots,\hat{t}_p,\ldots,t_k)$ (omitting $t_p$)
where $0\leq t_j<|m_j|$ for $j<i$ and others could be arbitrary
nonnegative integers.
\begin{itemize}
\item[(1)]If $t_{p+1}<|m_p|$, we define
$\mathbf{t}'_{p+1}=(t'_1,\ldots,\hat{t'}_{p+1},\ldots,t'_k)$ where
$t_j'=t_j$ for $j\neq p,p+1$ and $t_p'=t_{p+1}$. Obviously,
$\mathbf{t}'_{p+1} \in \daleth_{p+1}(m_1,\ldots,m_{p})$ and
$|\mathbf{t}'_{p}|=|\mathbf{t}'_{p+1}|$.

\item[(2)]If $t_{p+1}\geq|m_p|$, we define
$\mathbf{t}'_{p}=(t'_1,\ldots,\hat{t'}_{p},\ldots,t'_k)$ where
$t_j'=t_j$ for $j\neq p+1$ and $t'_{p+1}=t_{p+1}-|m_p|$, it's
clear to see that $\mathbf{t}'_{p} \in
\daleth_p(m_1,\ldots,m_{p-1})$ and
$|\mathbf{t}_p|=|\mathbf{t}'_{p}|+|m_p|.$
\end{itemize}
Furthermore, one can easily see that the above two  processes are
reservable. \qed

We now consider the coefficients of $\mathbf{x}^{\mathbf{m}_1}q^n$
in the series expansion of \eqref{rhand}. It's known that  only
the series expansion of the first term of \eqref{rhand} contains
the term $\mathbf{x}^{\mathbf{m}_1}q^n$, and the coefficients of
$\mathbf{x}^{\mathbf{m}_1}q^n$ are
\begin{align*}
&\sum_{\mathbf{t}_1}\mathcal{D}_1\left(m_1-\sum_{i=2}^k|m_i|-2|\mathbf{t}_1|-(k-1);n\right)\\
&=\sum_{j=0}^{+\infty}{j+k-2 \choose
k-2}N\left(\sum_{i=1}^k|m_i|+2j+k-1;n\right),
\end{align*}
where the equality follows from  the fact that $m_1\leq 0$,
$\mathcal{D}_1(-m;n)=N(m;n)$ and the number of solutions to
$t_2+t_3+\cdots+t_{k}=j$ in nonnegative integers is $j+k-2 \choose
k-2$.

Assume that the coefficients of $\mathbf{x}^{\mathbf{m}_p}q^n$ in
the series expansion of \eqref{rhand} equal to
\eqref{rghthandthm7}, we now show that the coefficients of
$\mathbf{x}^{\mathbf{m}_{p+1}}q^n$ are also equal to
\eqref{rghthandthm7}. Observe that $\mathbf{x}^{\mathbf{m}_p}q^n$
appears in the series expansions of the first $p$ terms of
\eqref{rhand} and the term $\mathbf{x}^{\mathbf{m}_{p+1}}q^n$
appears in the series expansions of the first $(p+1)$ terms.
Furthermore, the coefficients of $\mathbf{x}^{\mathbf{m}_p}q^n$
and $\mathbf{x}^{\mathbf{m}_{p+1}}q^n$ are the same in the series
expansions of the first $(p-1)$ terms. Therefore, if we verify the
sum of the coefficients of $\mathbf{x}^{\mathbf{m}_{p+1}}q^n$ in
the series expansions of the $p$th term and $(p+1)$th term of
\eqref{rhand} equal to the coefficients of
$\mathbf{x}^{\mathbf{m}_p}q^n$ in the series expansions of the
$p$th term, we could reach our conclusion by the induction
hypothesis.

By \eqref{coef-iterm}, it's known that the coefficients of
$\mathbf{x}^{\mathbf{m}_p}q^n$  (where $m_p\leq 0$) in the series
expansion of the $p$th term are
\begin{equation}\label{coef-pterm}
(-1)^{p-1}\sum_{\mathbf{t}_p}\mathcal{D}_1\left(\sum_{j=1}^{p-1}m_j-\sum_{j=p}^k|m_j|-2|\mathbf{t}_p|-(k-1);n\right),
\end{equation}
and the coefficients of $\mathbf{x}^{\mathbf{m}_{p+1}}q^n$ (where
$m_p\geq 1$ ) in the series expansions of the $p$th term and
$(p+1)$th term are
\begin{align*}
&(-1)^{p-1}\sum_{\mathbf{t}_p}\mathcal{D}_1\left(\sum_{j=1}^{p}m_j-\sum_{j=p+1}^k|m_j|-2|\mathbf{t}_p|-(k-1);n\right)\nonumber \\
&+(-1)^{p}\sum_{\mathbf{t}_{p+1}}\mathcal{D}_1\left(\sum_{j=1}^{p}m_j-\sum_{j=p+1}^k|m_j|-2|\mathbf{t}_{p+1}|-(k-1);n\right),
\end{align*}
which equal to \eqref{coef-pterm} by Lemma \ref{seq-lem-2}. Thus
we get our conclusion, and by the definition of
$R_k(x_1,x_2,\ldots,x_k;q)$ , we could recast \eqref{thm-7} as the
partition identity \eqref{main-e}.
 \qed

\section{$k$-marked odd Durfee symbols}
This section is devoted to solve the problems raised by Andrews
(\cite[Problems 6-9]{Andrews-07-a}) on $k$-marked odd Durfee
symbols. We begin this section by defining odd Durfee symbols
which resemble ordinary Durfee symbols with a modified subscript
and with odd numbers as entries.
\begin{defi}
An odd Durfee symbol of $n$  is a
 two-row array with subscript
\begin{equation}
(\alpha,\beta)_D=\left(\begin{array}{cccc}
\alpha_1&\alpha_2&\ldots&\alpha_s\\
\beta_1&\beta_2&\ldots&\beta_t
\end{array} \right)_{D}
\end{equation}
where $\alpha_i$ and $\beta_i$ are all odd numbers, $2D+1\geq
\alpha_1\geq \alpha_2\geq \cdots \geq \alpha_s>0$, $2D+1\geq
\beta_1 \geq \beta_2 \geq \cdots \geq \beta_t>0$, and
$n=\sum_{i=1}^s \alpha_i+\sum_{i=1}^t \beta_i+2D^2+2D+1$.
\end{defi}

The odd rank of an odd Durfee symbol is defined as the number of
parts of $\alpha$ minus the number of parts of $\beta$, let
$\mathcal{D}_1^0(m;n)$ denote the number of odd Durfee symbols of
$n$ with odd rank $m$, we then have

\begin{thm}\label{gf-r-odd}
The generating function for $\mathcal{D}_1^0(m;n)$ is given by
\begin{equation}
\sum_{n=0}^{+\infty}
\mathcal{D}_1^0(m;n)q^n=\frac{1}{(q^2;q^2)_{\infty}}\sum_{n=0}^{+\infty}(-1)^{n}q^{3n^2+3n+1+|m|(2n+1)}.
\end{equation}
\end{thm}
This result can easily follow by comparing the coefficients of
$z^m$ in  \eqref{gf-ds-odd} given by Andrews
\cite[(8.4)-(8.5)]{Andrews-07-a}:
\begin{align}
\sum_{n=1}^{+\infty}\sum_{m=-\infty}^{+\infty}\mathcal{D}_1^0(m;n)z^mq^n&=\sum_{n\geq
0}\frac{q^{2n(n+1)+1}}{(zq;q^2)_{n+1}(z^{-1}q;q^2)_{n+1}}\nonumber\\
&=\frac{1}{(q^2;q^2)_{\infty}}\sum_{n=-\infty}^{+\infty}\frac{(-1)^nq^{3n^2+3n+1}}{1-zq^{2n+1}},\label{gf-ds-odd}
\end{align}
where the first equality follows by direct combinatorial argument.
The second equality is given by \cite[p.66]{Watson-1936}.

Andrews also \cite{Andrews-07-a} defined the $k$th symmetrized odd
rank moment by
\begin{equation}\label{sym-oran-mom}
\eta^0_k(n)=\sum_{m=-\infty}^{{+\infty}}
{m+\lfloor\frac{k-1}{2}\rfloor \choose k}\mathcal{D}_1^0(m;n),
\end{equation}
and introduced $k$-marked odd Durfee symbols, whose definition is
almost identical to that of $k$-marked Durfee symbols (Definition
\ref{de-MDS}).

\begin{defi} A
$k$-marked odd Durfee symbol of $n$ is composed of $k$ pairs of
partitions into odd parts with the subscript,  which is defined as
 \[\eta^0=\left(\begin{array}{cccc}
\alpha^k,&\alpha^{k-1},&\ldots,&\alpha^1\\[2pt]
\beta^k,&\beta^{k-1},&\ldots,&\beta^1
\end{array}
\right)_D,\] where $\alpha^i$ {\rm(}resp. $\beta^i${\rm)} are all
partitions with odd parts and
$\sum_{i=1}^k(|\alpha^i|+|\beta^i|)+2D^2+2D+1=n$. Furthermore, the
partitions $\alpha^i$ and $\beta^i$ must satisfy almost the same
conditions with $k$-marked Durfee symbols  expect for the third
term in Definition \ref{de-MDS} where for the $k$th vector
$(\alpha^k, \beta^k)$ of $k$-marked odd Durfee symbol,
$\beta^k_1,\alpha^k_1 \leq 2D+1$.
\end{defi}
Following $k$-marked Durfee symbol, Andrews defined the $i$th odd
rank for  $k$-marked odd Durfee symbol.  For  a $k$-marked odd
Durfee symbol $\eta^0$,  we define $\rho_i(\eta^0)$, the $i$th odd
rank of $\eta^0$ by
\[\rho_i(\eta^0)=\left\{
\begin{array}{ll}
l(\alpha^i)-l(\beta^i)-1 \ \ &\text{ for }\ 1\leq i<k,\\[5pt]
l(\alpha^k)-l(\beta^k) \ \ &\text{ for }\ i=k.
\end{array}\right.
\]

Let $\mathcal{D}^0_k(m_1,m_2,\ldots,m_k;n)$ denote the number of
$k$-marked odd Durfee symbols of $n$ with $i$th odd rank equal to
$m_i$ and $\mathcal{D}^0_k(n)$ denote the number of $k$-marked odd
Durfee symbols of $n$. Define $R^0_k(x_1,x_2,\ldots,x_k;q)$ by
\begin{equation*}
R^0_k(x_1,x_2,\ldots,x_k;q)=\sum_{m_1,\ldots,m_k=-\infty}^{+\infty}\sum_{n=0}^{+\infty}
\mathcal{D}^0_k(m_1,m_2,\ldots,m_k;n)x_1^{m_1}x_2^{m_2}\cdots
x_k^{m_k}q^n.
\end{equation*}

Andrews deduced the following four identities on $k$-marked odd
Durfee symbols which are much similar with $k$-marked Durfee
symbols.

\begin{thm}[Corollary 27, Andrews \cite{Andrews-07-a}]\label{gf-km-odd}
\begin{align}\label{gf-km-odd-e}
R^0_k(x_1,x_2,\ldots,x_k;q)=\frac{1}{(q^2;q^2)_\infty}\sum_{n=0}^{+\infty}(-1)^{n}q^{3n^2+(2k+1)n+k}
\frac{1-q^{4n+2}}{\prod_{j=1}^k(1-x_jq^n)(1-\frac{q^n}{x_j})}.
\end{align}
\end{thm}

\begin{thm}[Corollary 28, Andrews \cite{Andrews-07-a}]\label{sym-odd}
$\mathcal{D}^0_k(m_1,m_2,\ldots,m_k;n)$ is symmetric in
$m_1,\,m_2,\ldots,m_k$.
\end{thm}
\begin{thm}[Corollary 29, Andrews \cite{Andrews-07-a}]\label{int-odd}
For $k\geq 1$,
\begin{equation}
\mathcal{D}^0_{k+1}(n)=\eta^0_{2k}(n).
\end{equation}

\end{thm}

\begin{thm}[Theorem 25, Andrews \cite{Andrews-07-a}]\label{thm-7-odd}
\begin{align}
R^0_k(x_1,x_2,\ldots,x_k;q)=\sum_{i=1}^k\frac{R^0_1(x_i;q)}{\prod_{\stackrel{j=1}{
j\neq i}}^k(x_i-x_j)(1-\frac{1}{x_ix_j})}.
\end{align}
\end{thm}

We  now give a brief expository of how to prove these four
conclusions combinatorially. First of all, it's straightforward to
see that the bijection $\Phi$ in Theorem \ref{r-kms-ds}, $\Psi$ in
Theorem \ref{r-km-kms} and the bijection $\Theta$ in Theorem
\ref{r-sym} on Section 2 are valid for $k$-marked odd Durfee
symbols, one then easily deduces the same result for $k$-marked
odd Durfee symbols as Theorem \ref{main}.
\begin{thm}\label{main-odd}
For $k\geq 2$, we have
\begin{equation}\label{main-odd-e}
\mathcal{D}^0_k(m_1,m_2,\ldots,m_k;n)=\sum_{j=0}^{+\infty}{j+k-2
\choose
k-2}\mathcal{D}_1^0\left(\sum_{i=1}^k|m_i|+2j+k-1;n\right).
\end{equation}
\end{thm}
Thus, Theorems \ref{sym-odd}, \ref{int-odd}, and \ref{thm-7-odd}
can be deduced from Theorem \ref{main-odd}  by the precisely same
progressions as Theorems \ref{sym-mds}, \ref{int}, and \ref{thm7}
on Section 3. To prove Theorem \ref{gf-km-odd}, it suffices to
prove that the right side hand of \eqref{gf-km-odd-e} is the
generating function for the summation on the right side of
\eqref{main-odd-e}, which can be easily derived by Theorem
\ref{gf-r-odd}, following the same progression as Theorem
\ref{gf-km}.

 \noindent{\bf Acknowledgments.} I would like to thank  Guoce Xin and Yue Zhou for helpful
discussions, and I am grateful to George E. Andrews for valuable
comments. This work was supported by the 973 Project, the PCSIRT
Project of the Ministry of Education, the Ministry of Science and
Technology, and the National Science Foundation of China.

\end{document}